\documentclass[11pt]{article}
\usepackage{graphicx}
\usepackage{amssymb,amsmath}
\usepackage{latexsym}
\usepackage{fullpage}
\usepackage{subfigure}

\usepackage{setspace}

\def\half{{\textstyle {\frac{1}{2}}}}

\def\R{\mathbb{R}}
\def\beq{\begin{equation}}
\def\eeq{\end{equation}}

\def\barr{\begin{array}}
\def\earr{\end{array}}
\def\ba*{\begin{eqnarray*}}
\def\ea*{\end{eqnarray*}}

\newcommand{\proofend}{\vrule height 6pt width 6pt depth 1pt}

\newtheorem{lemma}{Lemma}
\newtheorem{corollary}{Corollary}
\numberwithin{equation}{section}

\author{Lindsay B. H. May$^1$, Michael Shearer$^1$, Karen E. Daniels$^2$
\\ $^1$Department of Mathematics, North Carolina State University, Raleigh, NC 27695
\\ $^2$Department of Physics, North Carolina State University, Raleigh, NC 27695}

\title{Scalar conservation laws with nonconstant coefficients \\
with application to particle size segregation \\
in granular flow}

\date{July 29, 2009}

\begin{document}

\maketitle

\begin{abstract}
Granular materials will segregate by particle size when subjected to shear, as occurs, for example, in avalanches. The evolution of a bidisperse mixture of particles can be modeled by a nonlinear first order partial differential equation, provided the shear (or velocity) is a known function of position. While avalanche-driven shear is approximately uniform in depth, boundary-driven shear typically creates a shear band with a nonlinear velocity profile. In this paper, we measure a velocity profile from experimental data and solve initial value problems that mimic the segregation observed in the experiment, thereby verifying the value of the continuum model. To simplify the analysis, we consider only one-dimensional configurations, in which a layer of small particles is placed above a layer of large particles within an annular shear cell and is sheared for arbitrarily long times. We fit the measured velocity profile to both an exponential function of depth and a piecewise linear function which separates the shear band from the rest of the material. Each solution of the initial value problem is non-standard, involving curved characteristics in the exponential case, and a material interface with a jump in characteristic speed in the piecewise linear case.
\end{abstract}

\section{Introduction}

When set in motion through vibration or shear, granular materials have a strong tendency to segregate into bands  containing particles of similar size, shape, or density \cite{Ottino-2000-MSG}. In this paper, we focus on shear-induced segregation by size, which appears in a variety of configurations and applications, including avalanches \cite{Pouliquen-1999-SII}, rotating tumblers \cite{Metcalfe-1995-AMG}, and internal shear experiments \cite{Hill-2008-ISM}. 
 
Continuum models of avalanche flow have been derived using ideas from shallow water theory, in which a thin-layer approximation captures both the free surface shape and the underlying depth-averaged velocity \cite{savage-1989-mfm}. These models typically do not account for segregation. However, segregation in avalanching flows has been modeled by a mass transport equation alone, in which a roughly constant shear rate is specified \cite{graythornton1, savage}. Here, we adapt a mass transport  segregation model to situations where the shear rate is far from constant, reflecting the nonlinear dependence of  the velocity of particles on position, as is commonly the case for boundary-driven flows \cite{Midi-2004-DGF}. Through experiments on a mixture of two particle sizes within an annular shear cell, we measure the velocity profile and incorporate the resulting spatially-dependent shear rate into the constitutive law of the model.

The Gray-Thornton model \cite{graythornton1} for segregation by size in an avalanche containing two species of particles with similar density but different size, takes the form
 \beq\label{basic_pde}
 \varphi_t+u(z)\varphi_x+(w(\varphi,z)\varphi)_z=0.
 \eeq
In this partial differential equation (PDE), $\varphi(x,z,t)$ represents the concentration (fraction by volume) of smaller particles as a function of the distance down the avalanche $x$, distance $z$ above the base and time $t$.  The bulk flow is represented by the velocity $u(z)$ parallel to the base; the normal velocity $w(\varphi,z)$ of small particles is due to segregation dynamics. Both $u(z)$ and $w(\varphi,z)$ are assumed to be known functions whose functional forms are deduced as part of the model derivation. In avalanche flow, $u(z)$ is roughly linear near the surface \cite{Midi-2004-DGF}, so that the shear rate $|u'(z)|$ is close to constant, and $w(\varphi,z)= - k(1 - \varphi)$ for positive, constant $k$ proportional to the constant shear rate. Thus, in the Gray-Thornton model, the segregation rate is independent of depth. Note that $1-\varphi$ is the concentration of large particles, so that this form for $w(\varphi,z)$ may be considered to represent the availability of large voids created by relative motion of large particles.

In this paper, we are interested in the influence of non-uniform shear rate $|u'(z)|$, for which the segregation rate will be different at different depths. For example, if there is no shear, then there should be no tendency towards segregation, whereas a large shear rate should induce rapid segregation. Our model fits into the general framework of equation (\ref{basic_pde}), but the normal speed $w$ of small particles depends on $z$ through the depth-dependence of the shear rate $|u'(z)|$. In our model, we assume that $w$ is proportional to $|u'(z)|$. 
 In principle, $w$ could be any increasing function of shear rate. 

\begin{figure}
\begin{center}
\includegraphics[scale=.5]{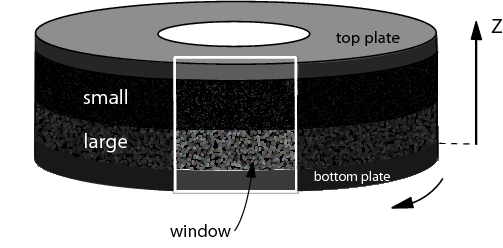}
\caption{Initial configuration of particle layers in the experimental annular Couette cell.} 
\label{exp1} 
\end{center}
\end{figure}

To simplify matters in both the experiment and the model, we begin with a bidisperse granular material in which the two sizes of particles (with the same density) are arranged in a one-dimensional configuration, as shown in Fig.~\ref{exp1}. It is reasonable to assume that in the subsequent evolution from this normally-graded configuration to an inverse-graded configuration, the concentration of each size of particle at each location depends only on depth and time. We explore two PDE models, both motivated by the structure of the velocity profile taken from experimental observations and differing only in the choice of depth-dependent shear rates chosen to approximate the experimental results.

The general form of the PDE (\ref{basic_pde}) we consider is the conservation law
 \beq\label{pde1}
 \varphi_t+(sa(z)f(\varphi))_z=0, \qquad 0<z<1, \quad t>0.
 \eeq
The convex flux $f:[0,1]\to \R$ satisfies $f(0)=f(1)=0$, corresponding to the physical limits of the dependent variable $\varphi$. The nonconstant coefficient $a(z)$ is the shear rate $|u'(z)|$; it is a monotonically decreasing function of the position variable $z$. The segregation rate parameter $s>0$ sets the time scale for the evolution of $\varphi$.

We will be concerned with initial boundary value problems, with initial data corresponding to the one-dimensional experimental configuration:
 \beq\label{ic1}
 \varphi(z,0)=\varphi_o(z)=\left\{\barr{ll}
 0, \quad &0<z<z_0,\\[10pt]
 1, & z_0<z<1,
 \earr
 \right.
 \eeq
 and boundary conditions 
 \beq\label{bcs1}
 \varphi(0,t)=1, \quad \varphi(1,t)=0.
 \eeq 
 
\bigskip

In the experiment, shown schematically in Fig.~\ref{exp1}, we place a layer of small glass spheres over a layer of larger spheres within the annular region between fixed rigid concentric cylinders. The aggregate is sheared by rotating the lower confining plate at fixed vertical position. An upper heavy confining plate is allowed to move vertically to accommodate changes in volume, due to both Reynolds dilatancy \cite{reynolds} and changes in packing density arising from the mixing/segregation process \cite{Golick-2009-MSR}.
The particles initially mix and then re-segregate through a process known as {\em kinetic sieving}: as the shearing proceeds, large particles roll and slide over one another, opening up gaps for the smaller particles to fall into. The small particles also act as levers for the large particles, which consequently tend to move vertically upwards, a process  sometimes called {\em squeeze expulsion}. Kinetic sieving was modeled in avalanche flow by Savage and Lun \cite{savage}, and subsequently by Gray and Thornton, using a different approach \cite{graythornton1}. In these models (which are valid in several space dimensions), the segregation rate is assumed to depend only on the concentration of small particles. This approximation is suitable for free-surface avalanches, where shearing is provided by the effect of gravity, a body force. However, in our experiments, shearing is instead provided by motion at the lower boundary, and this is transmitted through the granular material only by particle-particle contacts. The resulting shear rate drops off dramatically within a few layers of particles.

We model the depth-dependence of the shear rate $a(z)$ in two ways, suggested by the experimental data: 
\paragraph{Case~I: Piecewise constant shear rate:}
 \beq\label{case1}
 a(z)=\left\{\barr{ll}
 k_0, \quad &0<z<z_c,\\[10pt]
 k_1, & z_c<z<1,
 \earr
 \right.
 \eeq
with $k_0>k_1>0$.
\paragraph{Case~II: Smooth shear rate:} 
 \beq\label{case2}
 a(z)=a_0e^{-z/\lambda}, \qquad 0<z<1.
 \eeq

Case~I is based on the observation that there is a higher shear rate near the bottom plate, reflecting localization within a shear band. Modeling this higher rate as a constant is a coarse approximation to the experimental data. However, the split into two regimes, with a material interface at $z=z_c$ is justified by the data. Equation (\ref{pde1}) with a discontinuous function $a(z)$ does not fit into the existence theory of Kruzkov \cite{kruzkov}, and indeed, the issues of existence and uniqueness for this type of equation have been addressed in some generality only recently \cite{chen1}. In this case, characteristics are straight lines, on which $\varphi$ is constant, but both the characteristic speed and $\varphi$ experience a jump at $z=z_c$.

For smooth functions $a(z)$ (Case~II), the existence result of Kruzkov \cite{kruzkov} for initial value problems can be adapted to the initial boundary value problem, by extending $a(z)$ and the initial conditions beyond the boundary: $a(z)=a(0), \varphi(z,0)=1,  \ z<0; \ a(z)=a(1),  \varphi(z,0)=0, \ z>1.$  However, characteristics are curved, and moreover, $\varphi$ is not constant on characteristics. Consequently, although the structure of solutions can be characterized, the solutions cannot be found explicitly. The exponential form in Case~II is consistent with other studies of sheared granular materials \cite{Midi-2004-DGF}, and provides a close fit to our experimental data. We begin the analysis of Case~II by considering general functions $a(z)$ that are smooth, positive and decreasing, but it turns out that the choice of the exponential form is particularly useful for calculating explicit solutions.

Since our objective is to mimic the experiment, we restrict attention to initial conditions (\ref{ic1}) that reflect experimental conditions and for which we can analyze the solutions. We compute these solutions in detail using the structure of hyperbolic waves. Quantitative comparison of the theoretical solutions of this paper with experimentally-observed segregation is presented in \cite{May-2009-SDP}. 

\begin{figure}
\centerline{\includegraphics[width=0.95\textwidth]{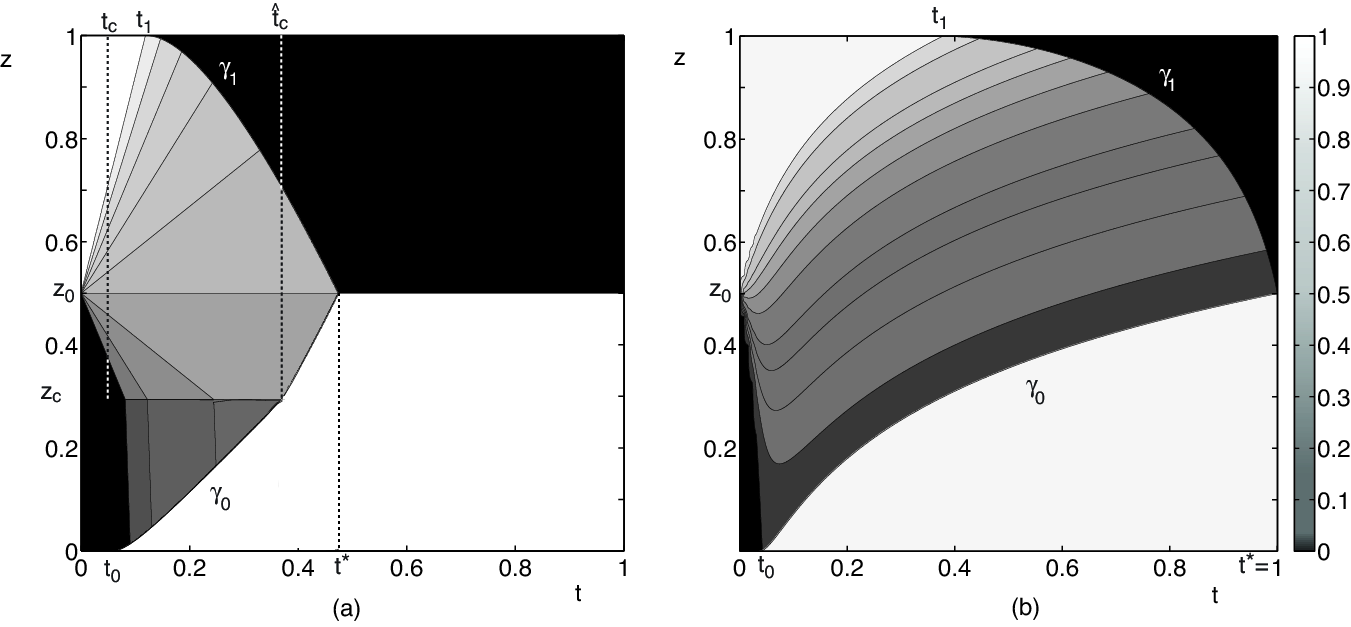}}
\caption{Solutions of (\ref{pde1}--\ref{bcs1}) in  (a) Case~I and  (b) Case~II with $s=13.60$ (which sets $t^* = 1$ in Case~II). }
\label{solutions} 
\end{figure}

The constants $z_c<z_0$ and $k_0>k_1>0$ in Case~I, and positive constants $a_0, \lambda$ in Case~II, are determined from an experimentally-measured velocity profile, and an overall segregation rate constant $s$ sets the time scale. Solutions based on the experimentally determined constants are shown in Fig.~\ref{solutions}, in which $s$ is chosen to make the final time $t^*=1$ in Case~II. The solutions involve a rarefaction wave, centered at $(z,t)=(z_0,0)$, in which $\varphi(z,t)$ varies continuously between $\varphi=0$ and $\varphi=1$. As the leading edge reaches $z=1$ (at time $t = t_1$), a shock wave is reflected, with a layer of large particles ($\varphi=0$) growing behind it. Similarly, as the trailing edge hits the boundary $z=0$ (at time $t = t_0$), a layer of small particles ($\varphi=1$) develops behind the reflected shock. The two shocks eventually meet at time $t^*$, at which time the solution becomes a stable stationary shock representing a layer of large particles above a layer of small particles separated at $z=z^*=1-z_0$. 

In \S\ref{sec_exp1}, we describe the annular shear cell experiment and explain how we determine the model parameters for the two chosen cases (\ref{case1}), (\ref{case2}). In \S\ref{IBVP}, we construct solutions in each of the two cases. Interestingly, in Case~I, the time to full segregation is independent of the shear rate $k_0$ within the shear band. We conclude with a discussion in \S\ref{discuss}.

\section{Experimental Results}\label{sec_exp1} 

The experimental configuration is an annular Couette cell (see Fig.~\ref{exp1}) with channel width 3.8~cm bounded by concentric aluminum cylinders with inner and outer radii 25.5~cm and 29.3~cm, respectively. The rotating bottom plate and an upper confining plate each have rubberized surfaces to enhance friction with the particles. A motor drives the bottom plate at a constant rotation rate of approximately 3 revolutions per minute. The cell is filled with a 2~kg layer of glass spheres (diameter 3~mm), placed above a 2~kg layer of larger glass spheres (diameter 6~mm). The fill height is approximately 4.1~cm, and a heavy top plate confines the particles but is free to move vertically to accommodate changes in volume as the aggregate dilates, mixes and segregates. Further experimental details are available in \cite{Golick-2009-MSR, May-2009-SDP}.

The apparatus has a window in the outer wall, permitting us to track particle positions over time with a high speed (450 Hz) digital camera. The camera collects digital images at discrete intervals throughout the duration of the experiment, allowing us to compare particle velocities at different stages of the experiment. In each  image, we locate the center of each particle, distinguishing large from small, and record the positions of individual particles. Through an automated process, we identify the same particle in successive frames, generating a list of the horizontal and vertical coordinates of each particle at a sequence of times. We refer to this list as a single-particle trajectory. 

\begin{figure}
\centering
\subfigure[]{
\includegraphics[width=0.47\textwidth]{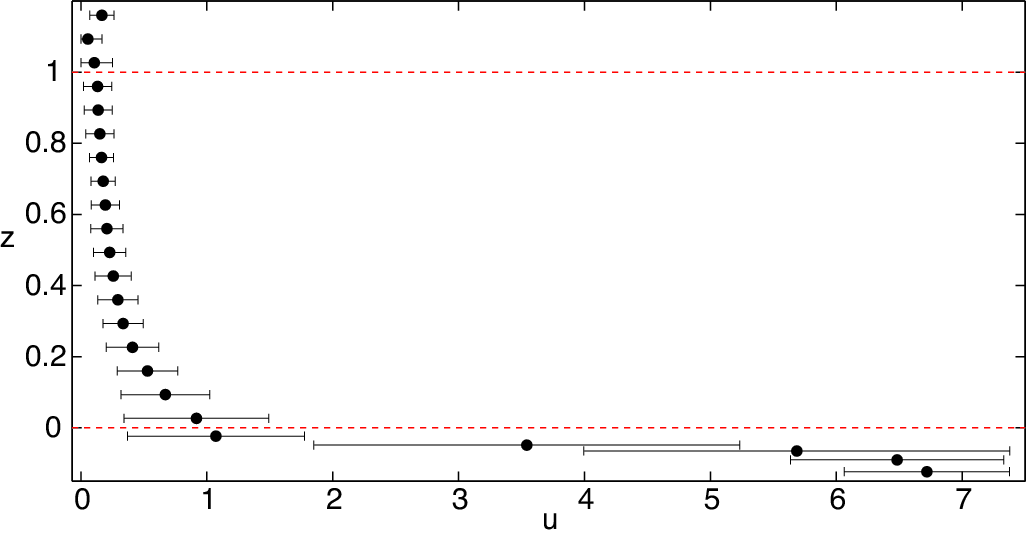}
\label{VP:full}
}
\subfigure[]{
\includegraphics[width=0.47\textwidth]{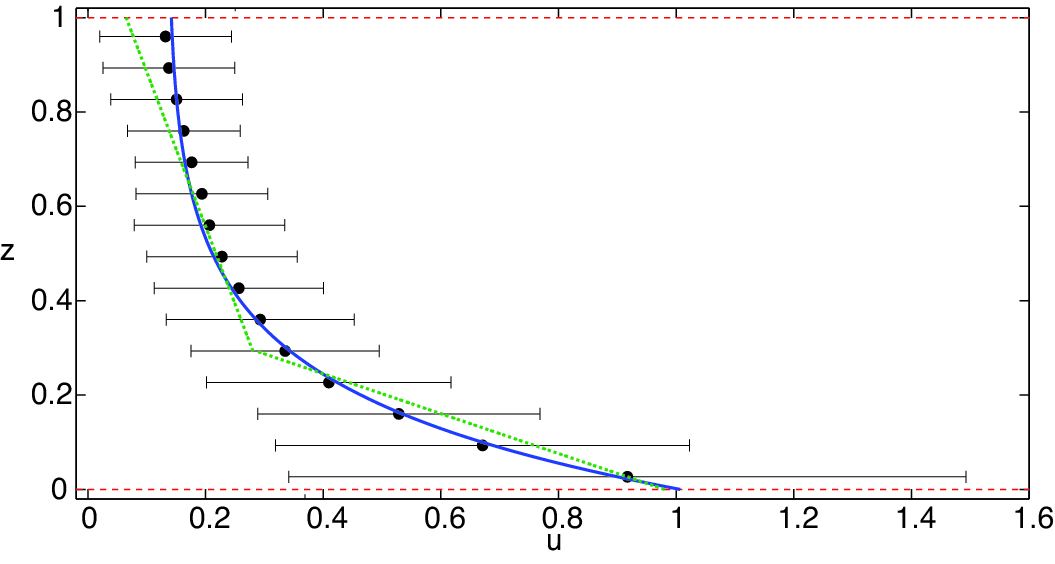}
\label{VP:zoom} 
}
\caption{Measured velocity profile $u_i$ ($\bullet$) for \subref{VP:full} full cell height, with boundary layers above and below the dashed horizontal lines and \subref{VP:zoom} within the region $z = [0,1]$. The dotted line is the fit to Case~I; the solid line is the fit to Case~II, as described in \S\ref{findingparameters}. The velocity profile is scaled so that $u(0) = 1$. }
\label{VP}
\end{figure}

\begin{figure}
\centering
\subfigure[]{
\includegraphics[width=0.47\textwidth]{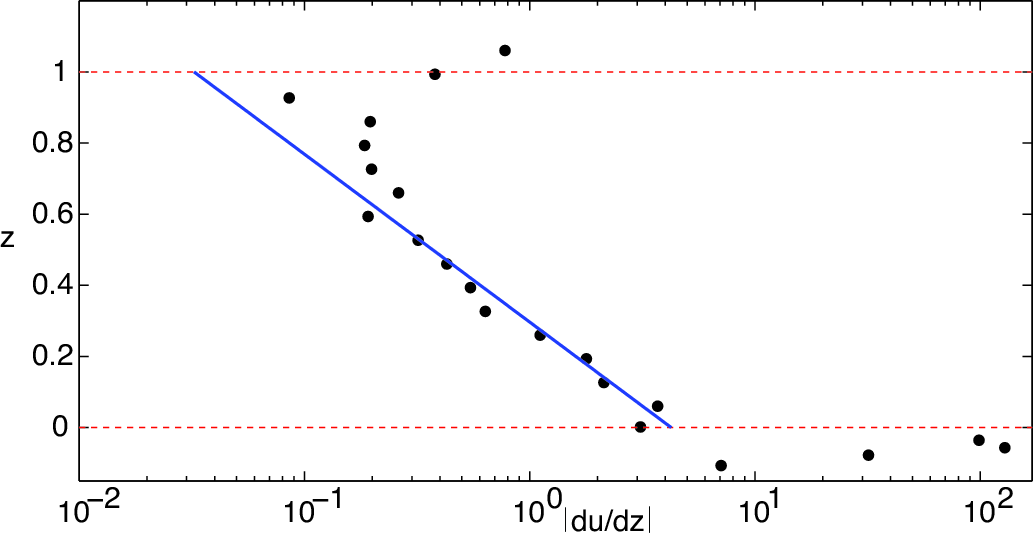}
\label{SR:semilog}
}
\subfigure[]{
\includegraphics[width=0.47\textwidth]{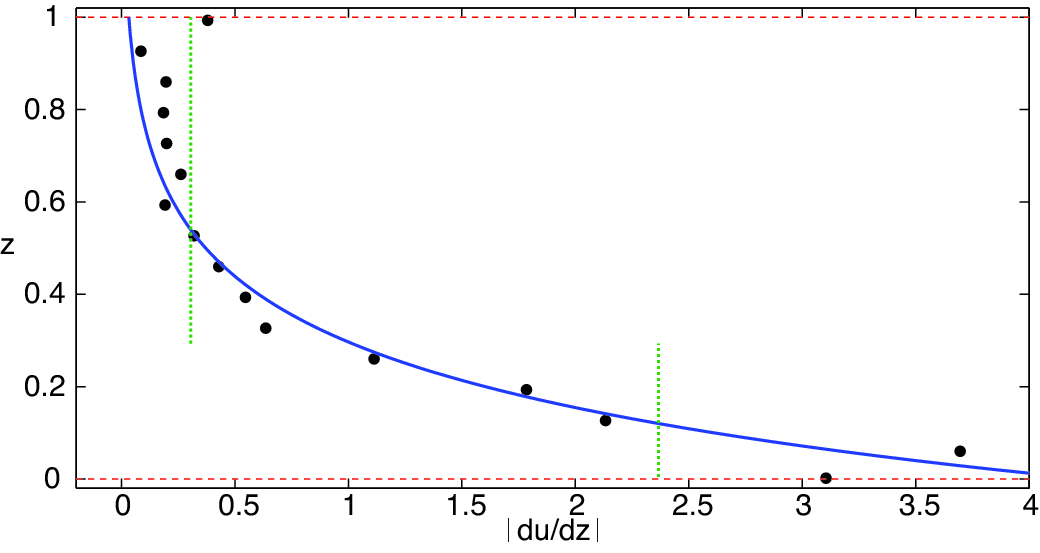}
\label{SR:fits} 
}
\caption{Shear rate $|du/dz|$ ($\bullet$) for \subref{SR:semilog} full cell height, with boundary layers above and below the dashed horizontal lines and \subref{SR:fits} within the region $z = [0,1]$. Vertical dotted lines are the fit to Case~I; solid line is the fit to Case~II, as described in \S\ref{findingparameters}.}
\label{SR}
\end{figure}

For each single-particle trajectory, we calculate the instantaneous horizontal velocity of the particle as follows. First, the vertical dimension of the sample is divided into twenty-three bins centered at positions $z=z_i, \, i=1, \ldots, 23$. Each trajectory is assigned to a bin $z_i$ based on the average vertical position of the particle. Using a moving interval of duration $\Delta t$, we determine the instantaneous velocity by fitting a linear function to the horizontal coordinates within $\Delta t$. For each bin, we choose an appropriate, speed-dependent, value for $\Delta t$, which varies from approximately $0.1$ seconds near the bottom plate to approximately $0.4$ seconds near the top plate. This process yields a range of velocities observed for an ensemble of different particles at different times. We fit a parabola to the peak of the probability distribution within each bin to calculate a velocity $u_i$ representing the horizontal speed of particles in bin $i$.

Fig.~\ref{VP}(a) shows the velocities $u_i, \, i=1, \ldots, 23$, which we refer to as the {\em measured velocity profile.} In the figure, we have normalized $z$ to $[0, 1]$ over the region of interest (described in \S\ref{findingparameters}), and scaled the velocity so that $u=1$ at $z=0$. The error bar through each point $(u_i, z_i)$ is the width of the parabola at a height one half of the maximum height, to give a sense of the distribution of observed values in each bin. In total, the measured velocity profile is based on processing particle positions from approximately $7\times 10^5$ images. In \S\ref{findingparameters}, we use the measured $u_i$ to generate appropriate parameters for the shear rate $a(z)$ for use in the model. Further details concerning the collection and processing of the experimental data are described in \cite{May-2009-SDP}.

During the processing of the data to generate the measured velocity profile, we established two properties which are crucial for the continuum model: 
\begin{itemize}
\item Velocities are similar for both large and small particles; calculating $u_i$ separately for large and small particles gives negligible differences. 
\item Velocities reach steady-state after a short initial transient of approximately $ 0.05 \, t^*_\mathrm{exp}$, where $t^*_\mathrm{exp} = 700$~seconds is the duration of the experiment. This observation justifies the use of a time-independent velocity profile $u(z)$ in the model.
\end{itemize}

\subsection{Determining the Shear Rate Profiles} \label{findingparameters}

To obtain the position-dependence of the shear rate from the measured velocity profile, we first take finite differences of $u_i$ between adjacent layers $z_i$:
\beq\label{shearrates}
u'(z_{i+\half})\approx \frac{u_{i+1}-u_{i}}{\Delta z_i}, \quad \Delta z_i= z_{i+1}-z_i, \quad z_{i+\half}=\half(z_i+z_{i+1}).
\eeq
The resulting shear rates are shown as solid points in Fig.~\ref{SR}. 

In Fig.~\ref{SR}\subref{SR:semilog}, we observe that the shear rates naturally fall into three sections, marked by the horizontal dashed lines in both Fig.~\ref{VP} and Fig.~\ref{SR}. The uppermost ($z>1$) and lowermost ($z<0$) regions are the boundary layers.  When the height of the sample is measured in real units, the width of each boundary layer is equivalent to one large particle diameter or two small particle diameters.  We employ a linear transformation to ensure that  $z=0$ and $z=1$ correspond to the top of the lower boundary layer and the bottom of the upper boundary layer, respectively. We limit our modeling to $z \in [0,1]$ since we are interested in the bulk behavior of the system.  We also normalize the velocity so that $u(0)=1.$  Since there is no data point at $z=0,$ we calculate the line containing the points $(u_5, z_5)$ and $(u_6,z_6),$ which span $z=0$, and use the velocity value associated with $z=0$ on that line to normalize the velocity data 
(see Fig.~\ref{VP}\subref{VP:full}). Note that the velocity of the bottom plate sets an overall timescale that is necessary to make a full comparison between predictions of the model and the observed segregation in the experiment. However, in this paper we consider only a comparison between the theoretical predictions of Case~I and Case~II, using the experiment solely to provide physically realistic shear rate parameters.

In Fig.~\ref{SR}\subref{SR:fits} we observe that the shear rates can be split into a low-shear region and a high-shear region. The division occurs at $z_c$, which we take to be located midway between   two adjacent $z_{i+\half}$ points: 
 $z_c=z_{10}= 0.29$. 
 To determine shear rate parameters $k_0$ and $k_1$, we average the shear rates in each of the two regions. This yields
$k_0= 2.4 \pm 0.4$  
for $0<z<z_c$ and
$k_1= 0.31 \pm 0.05$  
for $z_c<z<1$, which are both shown as vertical dotted lines in Fig.~\ref{SR}\subref{SR:fits}. For comparison, we can use $(k_0, k_1, z_c)$ to generate the corresponding piecewise linear fit to the measured velocity profile; this is shown by the dotted lines in Fig.~\ref{VP}\subref{VP:zoom}. These three parameter values are used with the constructions of \S\ref{IBVP} to generate the solution in Case~I shown in Fig.~\ref{solutions}(a).

The measured velocity profile in the region $0\leq z\leq 1$ is also well-described by an exponential function $u(z)=be^{-z/\lambda}+c$, as shown in Fig.~\ref{VP}\subref{VP:zoom}. A least-squares fit provides model parameters
 $\lambda = 0.22\pm 0.01$ and $b = 0.82\pm 0.05$
 and $c = 0.14 \pm 0.06.$ 
The resulting shear rate $|u'(z)| = \frac{b}{\lambda} e^{-z/\lambda}$ is plotted as a straight line in the semi-logarithmic plot Fig.~\ref{SR}\subref{SR:semilog} and as a curved line in Fig.~\ref{SR}\subref{SR:fits}; these figures verify that the procedure for determining the exponential shear rate from the measured velocity profile also provides a good fit to the shear rates. The parameter values $\lambda$ and $b$ are used with the constructions of \S\ref{IBVP} to generate the solution in Case~II shown in Fig.~\ref{solutions}(b).

To summarize, we have determined parameter values from the experiment for shear rates in Case~I and Case~II.
The specific values we use to generate the solutions shown in Figure~\ref{solutions} are:
\beq\label{values}
\barr{ll}
\mbox{Case~I:}\quad &z_0=0.5, \quad z_c=0.29, \quad k_0=2.4, \quad  k_1=0.31. \\[10pt]
\mbox{Case~II:} \quad & 
z_0=0.5, \quad \lambda=0.22, \quad b=0.82,\quad  a_0=b/\lambda=3.7.
\earr
\eeq

\section{Initial Boundary Value Problems} \label{IBVP}

In this section, we derive solutions of the initial boundary value problem (\ref{pde1}--\ref{bcs1}) in Cases~I and II. Since the segregation rate parameter $s>0$ simply affects the time scale, we first set $s=1$, and later normalize the time scale by choosing the value for $s$ which provides $t^*=1$ in Case~II. We begin with a treatment of characteristics and shocks, focusing on differences from standard constructions.

\subsection{Characteristics and Shocks}\label{chars1} 

Characteristics reduce the construction of continuous solutions of scalar first order PDEs to solving ordinary differential equations. For equation (\ref{pde1}) with $s=1$, characteristics are curves $z=z(t)$ and $\varphi=\varphi(t)$ given by 
\beq\label{chars1_1}
\frac{dz}{dt}=a(z)f'(\varphi); \quad \frac{d\varphi}{dt}=-a'(z)f(\varphi).
\eeq
Thus, $a(z)f(\varphi)$ is conserved along characteristics:
\beq
\label{invariant}
a(z)f(\varphi)= \mbox{constant}.
\eeq

Along characteristics in Case~I, in which $a(z)$ is piecewise constant, $z(t)$ is piecewise linear with a jump in slope across $z=z_c$ and $\varphi$ is piecewise constant with a jump across $z=z_c$. In Case~II, the characteristics are smooth curves: since $a'(z) \neq 0$, the only characteristics which are straight lines are those with $\varphi=0$ or $\varphi=1$. All other characteristics are not straight, and moreover, $\varphi$ is not constant along them. This is in agreement with the observation that $\varphi=0$ and $\varphi=1$ are the only constant solutions of the PDE in Case~II.

\begin{figure}
\begin{center}
\includegraphics[scale=.8]{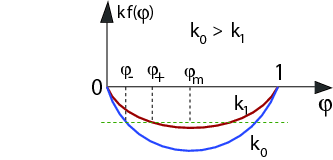}
\caption{The flux balance condition (\ref{flux0}) in Case~I.} 
\label{twoks} 
\end{center}
\end{figure}

Shock waves satisfy the Rankine-Hugoniot condition, in which the speed of the shock is related to the flux across it. Specifically, if the shock is $z=\gamma(t)$, and $\varphi_\pm(t)=\varphi(\gamma(t)\pm,t), a_\pm(t)=a(\gamma(t)\pm)$ are the one sided limits, then
\beq\label{rh1}
\frac{d\gamma}{dt}=\frac{a_+(t)f(\varphi_+(t))-a_-(t)f(\varphi_-(t))}{\varphi_+(t)-\varphi_-(t)}.
\eeq 
This formulation is consistent with the interpretation of the interface $z=z_c$ in Case~I as a stationary shock, for which $\gamma(t)=z_c$, and across which the fluxes balance:
\beq\label{flux0}
k_0f(\varphi_-(t))=k_1f(\varphi_+(t)).
\eeq
The flux balance (\ref{flux0}) is also consistent with the structure (\ref{invariant}) of characteristics. In Fig.~\ref{twoks} we show typical fluxes in Case~I with values of $\varphi_\pm$ representative of the solution of our specific initial value problem. In the figure, $f(\varphi)$ has a minimum at $\varphi=\varphi_m$. In our case, $\varphi_+<\varphi_m$ is known, and $\varphi_-<\varphi_+$ is then determined from (\ref{flux0}). However, if $\varphi_-<\varphi_m$ were given and $k_0f(\varphi_-)<k_1f(\varphi_m)$, then there would be no value of $\varphi_+$ satisfying (\ref{flux0}). Consequently, the solution of the initial value problem would be rather different, with a shock wave reflected from the interface $z=z_c$.

\subsection{Case~I} 

In this subsection, we solve the initial boundary value problem (\ref{pde1}--\ref{bcs1}) in Case~I, in which $a(z) $ is given by (\ref{case1}). For the solution, it is crucial that $k_0>k_1$, which is the physically-meaningful relationship for granular shear bands. In addition, since $z_0=\half$ in the experiment, we also assume $z^*=1-z_0>z_c$. In Fig.~\ref{solutions}(a), we show the solution with values of $k_0$ and $k_1$ calculated from the measured velocity profile in \S\ref{findingparameters}.

We construct the solution $\varphi(z,t)$ in several steps, corresponding to the different features in Fig.~\ref{solutions}(a). For small $t>0$, the solution consists of a single rarefaction wave centered at $(z,t)=(z_0,0)$. The rarefaction reaches the material interface $z=z_c$ at a time $t=t_c$, and is transmitted through through the interface as a simple wave, in general not centered. The simple wave first reaches the boundary $z=0$ at a time $t=t_0$, and the rarefaction wave first reaches the boundary $z=1$ at a time $t=t_1$. The simple wave and rarefaction are reflected from the boundaries as shock waves $z=\gamma_j(t), \, t>t_j, \, j=0,1$.  The shock $z=\gamma_0(t)$ crosses the interface $z=z_c$ at a time $t=\hat{t}_c$, and meets the shock $z=\gamma_1(t)$ at a time $t^*$. For $t>t^*$, 
 the solution is the piecewise constant function 
$$
\varphi(z,t)=\left\{\barr{ll} 1, \quad & z<z^*\\[10pt]
0&z>z^*,
\earr
\right.
$$ 
where $z^*=1-z_0$, as expected from conservation of the total mass (or volume) of small particles. Note that if $z^*<z_c$, then the descending shock $z=\gamma_1(t)$ reaches the interface $z=z_c$ before the rising shock $z=\gamma_0(t)$, and is transmitted through the interface; apart from this difference, the solution is the same.

To simplify some of the construction, and carry explicit calculations as far as possible, we restrict attention to the case $f(\varphi)=\varphi(\varphi-1)$. Then, the centered rarefaction is given explicitly by
\beq \label{phiztMain}
\varphi=\varphi_1(z,t) = \frac{1}{2}\left(\frac{z-z_0}{k_1t}+1\right),\quad z_c<z<1, \, t>0. 
\eeq
It reaches the boundary $z=1$ at time $t_1=(1-z_0)/k_1$, since the first characteristic to reach this boundary carries $\varphi=1$. The rarefaction is reflected as a shock $z=\gamma_1(t)$ satisfying the jump condition (\ref{rh1}). Since the shock has the centered rarefaction on one side, and $\varphi=0$ on the other, it is determined from the initial value problem
\beq 
\frac{d\gamma_1}{dt}=k_1(\varphi-1)=\frac{\gamma_1}{2t}-\frac{z_0}{2t}-\frac{k_1}{2}, \quad \gamma_1(t_1)=1.
\eeq
Thus,
\beq
\gamma_1(t)=z_0 - k_1t + 2\sqrt{(1-z_0)k_1t}.
\eeq

Similarly, the characteristic with $\varphi=0$ reaches $z=z_c$ at time $t_c = (z_0 - z_c)/k_1$. Along the line $z=z_c$ and $\varphi$ takes values 
\beq\label{phiplus}
\varphi_+(t)=\frac 12\left(1-\frac{t_c}{t}\right).
\eeq
The line $z=z_c$ in the $(z,t)$-plane behaves as a stationary shock as far as the weak solution is concerned. Consequently, $\varphi$ jumps from $\varphi_+$ to a value 
$\varphi_-(t)=\varphi(z_c-,t), \, t>t_c$, while keeping the flux continuous; the jump condition (\ref{flux0}) is 
\beq \label{flux1}
k_0\varphi_-(\varphi_--1) = k_1\varphi_+(\varphi_+-1).
\eeq
Solving this quadratic equation for $\varphi_-\in [0,\half)$, we find
 \beq\label{phiminus}
 \varphi_-(t) = \frac{1}{2}\left( (1-\sqrt{4\frac{k_1}{k_0}\varphi_+(\varphi_+-1)} \right),
\eeq
with $\varphi_+=\varphi_+(t)$ given by (\ref{phiplus}).

Next, we construct the simple wave that emanates from the line $z=z_c$. The construction involves a family of straight line characteristics parameterized by $\tau\geq t_c:$ 
\beq\label{zinR1}
z=k_0(2\varphi-1)(t-\tau)+z_c, \quad  t>\tau.
\eeq
 On each characteristic, $\varphi=\varphi_-(\tau)$ is constant. Thus,  
equations (\ref{phiplus}), (\ref{phiminus}), (\ref{zinR1}) define $\varphi(z,t)$ implicitly in the simple wave.

 In order to calculate the shock wave $z=\gamma_0(t)$ that reflects from the boundary $z=0$, we need to be able to calculate $\varphi(z,t)$ in the simple wave. Apart from the outermost characteristic
\beq\label{char_0}
 z=-k_0(t-t_c)+z_c,
 \eeq
on which $\varphi=0$ is constant, we find $\varphi(z,t)$ numerically by solving a quartic equation, derived as follows. 

The function $\varphi_-(\tau)$ is defined by (\ref{phiplus}), (\ref{phiminus}). We can write the inverse of this function, obtaining $\tau=\tau(\varphi)$ :
 \beq\label{tau}
 \tau(\varphi) = \left(\sqrt{\frac{k_1}{k_0(2\varphi-1)^2-k_0+k_1}}\right)t_c. 
 \eeq
Substituting into equation (\ref{zinR1}), we have an equation defining $\varphi$ as a function of $z$ and $t$. 
Let $\psi = 2\varphi-1$, $\alpha = 1-\frac{k_1}{k_0} >0$, and $\beta = \left(\sqrt{\frac{k_1}{k_0}}\right)t_c > 0$. Then in the new parameters and variables, (\ref{zinR1}), (\ref{tau}) become
\begin{eqnarray*}
z = k_0\psi\left(t-\frac{\beta}{\sqrt{\psi^2-\alpha}}\right)+z_c.
\end{eqnarray*}
Rearranging and expanding, we find that we have a quartic equation for $\psi$:
\begin{eqnarray}
g(\psi;z,t) \equiv A\psi^4+B\psi^3+C\psi^2+D\psi+E=0,
\label{gofpsi}
\end{eqnarray}
with coefficients depending on $z,t$ given by
\beq
A \equiv k_0^2t^2, \,
B \equiv 2k_0t(z_c-z), \,
C \equiv (z_c-z)^2-k_0^2(\alpha t^2+\beta^2), \,
D \equiv -2k_0t(z_c-z)\alpha, \,
E \equiv -\alpha(z_c-z)^2.
\eeq
For $(z,t)$ in the simple wave, we seek to solve equation (\ref{gofpsi}) for $\psi\in (-1,-\sqrt{\alpha})$, corresponding to $0<\varphi<\half(1-\sqrt{\alpha})$. The solution can then be used to find the shock wave $z=\gamma_0(t)$.

\bigskip

\begin{lemma}
For $(z,t)$ in the simple wave, $g(-1;z,t)\geq 0>g(-\sqrt{\alpha};z,t)$, with $g(-1;z,t)=0$ only on the characteristic (\ref{char_0}).
\end{lemma}
{\bf Proof:}
It is straightforward to check $g(-1;z,t)=0$ on the characteristic (\ref{char_0}), so we suppose that $(z,t)$ lies above that characteristic in the $(z,t)$ plane.
First we will show $g(-1;z,t) = A - B + C - D + E > 0$. Substituting in the values for the coefficients and simplifying, we find 
\begin{eqnarray*}
g(-1;z,t) = k_0k_1(t^2-t_c^2)+2k_1t(z-z_c)+\frac{k_1}{k_0}(z_c-z)^2.
\end{eqnarray*}
Then we substitute for $z$ using equation (\ref{zinR1}) and simplify, concluding that 
\begin{eqnarray*}
g(-1;z,t) = k_0k_1\left(\tau^2-t_c^2+4\varphi\tau(t-\tau)+4\varphi^2(t-\tau)\right). 
\end{eqnarray*}
Along the characteristic (\ref{zinR1}) in the simple wave, we have $t>\tau$, and $\tau>t_c$ except along the straight characteristic along which $\varphi=0$. Therefore, $g(-1;z,t)>0$. 

 At $\psi=-\sqrt{\alpha}$, a similar calculation yields 
\begin{eqnarray*}
g\left(-\sqrt{\alpha};z,t\right) = -k_1\left(k_0-k_1\right)t_c^2<0,
\end{eqnarray*}
since $k_0>k_1>0$. This completes the proof.  \proofend

\bigskip

\begin{corollary}
For $(z,t)$ in the simple wave, $g(\psi;z,t)=0$ has a solution in the interval $[-1,-\sqrt{\alpha})$, and a positive solution.
\end{corollary}
{\bf Proof:} \ The Lemma establishes the solution in $[-1,-\sqrt{\alpha})$. Since the constant $E$ in (\ref{gofpsi}) is negative, the product of the four solutions of $g(\psi;z,t)=0$ is negative. Thus, whether the polynomial has all real roots, or two real and two complex conjugate roots, at least one of the roots must be positive, since we already have established a negative root.
\proofend

\bigskip

The leading edge of the simple wave is the characteristic (\ref{char_0}) on which $\varphi=0$. 
 It reaches the boundary $z=0$ at time $t=t_0$ given by $t_0=t_c+ \frac{z_c}{k_0}$. 
From the point $(z,t)=(0, t_0)$, a shock $z=\gamma_0(t)$ emerges from the boundary. Behind the shock is a layer of small particles, with $\varphi=1$. Consequently, from the Rankine-Hugoniot condition (\ref{rh1}), the reflected shock satisfies
\beq\label{shock1}
\frac{d\gamma_0}{dt} = k_0\varphi(\gamma_0,t), \quad \gamma_0(t_0) =0,
\eeq
where $\varphi(\gamma_0,t)$ is the value of $\varphi$ in the simple wave at the shock. 

Equation (\ref{shock1}) is solved numerically, since we do not have a closed formula for the simple wave. The Corollary shows that $\varphi(\gamma_0(t),t)$ can be determined by solving equation (\ref{gofpsi}). To solve equation (\ref{shock1}), we therefore use the {\scshape Matlab} function {\tt roots} in conjunction with the {\scshape Matlab} routine {\tt ode45}, employing the values (\ref{values})
determined from the experimental data in \S\ref{findingparameters}. At each call of {\tt roots}, we verify that $g(\psi;z=\gamma_0(t),t)$ has two complex roots, thereby checking that we have found the only relevant value of $\varphi$ in the simple wave.

\begin{figure}
\begin{center}
\includegraphics[scale=.7]{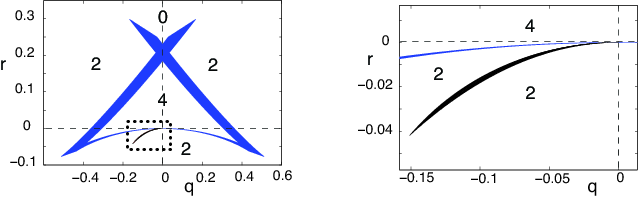}
\caption{Projection of the swallowtail and coefficients of $g(\psi;z,t)$ in the simple wave in Case~I, showing the number of real solutions of equation (\ref{swallow1}). The figure on the right is a magnification of the dotted rectangle in the left figure. } 
\label{swallow4} 
\end{center}
\end{figure}

As a further check, we compare the coefficients in equation (\ref{gofpsi}) with an established criterion for the existence of just two real roots. To do so, we place the quartic equation into a normal form
\beq\label{swallow1}
x^4+px^2+qx+r=0,
\eeq
by dividing (\ref{gofpsi}) by the coefficient $A$, and letting $x=\psi+\frac{B}{4A}$. The coefficients $p,q,r$ are then functions $\hat{p},\hat{q},\hat{r}$ of $(z,t)$, in addition to the parameters $k_0,k_1,z_c,z_0$. Equation (\ref{swallow1}) has coincident roots on the swallowtail surface $\cal{S}$ generated by eliminating $x$ from equation (\ref{swallow1}) and the equation
\beq\label{swallow2}
4x^3+2px+q=0.
\eeq
A convenient parametrization of $\cal{S}$ is obtained by expressing $(p,q,r)$ in terms of $p$ and $x$ :
\beq\label{swallow3}
\left\{
\barr{rcl}
q&=&-2px-4x^3\\
r&=&-qx-px^2-x^4=px^2+3x^4.
\earr
\right. 
\eeq
For the parameter values (\ref{values}), we easily verify that $\hat{p}(z,t)<0$ for $(z,t)$ in the simple wave 
$0<z<z_c, t>t_c$ and $(\hat{q},\hat{r})(z_c,t)\equiv 0$. Moreover, the surface $\hat{\cal{S}}=\left\{(\hat{p},\hat{q},\hat{r})(z,t): 0<z<z_c, t>t_c\right\}$ lies below the swallowtail $\cal{S}$. This region in $(p,q,r)$-space corresponds to coefficient values for which (\ref{gofpsi}) has exactly two real roots. In Fig.~\ref{swallow4}, we show the projection of the swallowtail onto the $(q,r)$ plane for values of $p$ including the range of $\hat{p}$. We superimpose the corresponding projection of $\hat{\cal{S}}$.

\ \\

Returning to the structure of the solution of the initial boundary value problem, we observe that, although the shock is initially tangent to the $t$-axis, it then immediately has positive speed,  since $\varphi(z,t)>0$ for $t>t_0, z\geq 0$ in the simple wave.
Consequently, the shock $z=\gamma_0(t)$ reaches $z=z_c$ at a finite time $t=\hat{t}_c.$ Since $\varphi_\pm=1$ satisfies the compatibility condition (\ref{flux1}), the shock $z=\gamma_0(t)$ is simply transmitted through the interface but now satisfying the ODE $\gamma_0'(t) = k_1\varphi(\gamma_0,t)$, with initial condition $\gamma_0(\hat{t}_c)=z_c$, and $\varphi=\varphi_1(z=\gamma_0,t)$ given by the centered rarefaction wave (\ref{phiztMain}). 
Solving the initial value problem, we find an explicit formula for the solution
\beq
 \gamma_{0}(t) =z_0+ k_1t+\sqrt{\frac{t}{\hat{t}_c}}(z_c-z_0-k_1\hat{t}_c), \quad t>\hat{t}_c.
\eeq
 
 To determine the time $t^*$ at which shocks $\gamma_{0}$ and $\gamma_1$ meet, we first note that by mass conservation, they meet at the location  $z=1-z_0$. Then $t^*$ can be determined from the equation $\gamma_1(t^*)=1-z_0,$ resulting in the expression
\beq\label{tstar}
t^*=\frac{1}{k_1}\bigg(1+2\sqrt{z_0(1-z_0)}\bigg).
\eeq
But then $\gamma_0(t^*)=1-z_0$ becomes an equation for $\hat{t}_c$, with the result that $\hat{t}_c$ is independent of $k_0,$ and can be calculated explicitly, without resorting to the numerical values of the shock $z=\gamma_0(t)$ as it approaches $z=z_c$ from below:
\beq\label{thatc}
\hat{t}_c=\frac{1}{k_1}\bigg(\sqrt{z_0}+\sqrt{z_c}\bigg)^2.
\eeq

In Fig.~\ref{solutions}(a), we show $\hat{t}_c$ and $t^*$ normalized by the segregation rate constant $s=13.60$. In the figure, we use the parameter values (\ref{values}).  Then  (\ref{tstar}), (\ref{thatc}) give
$$
\mbox{Case~I:} \qquad \hat{t}_c= 0.37; \quad t^* =0.47,
$$
in agreement with the simulation.

It is remarkable that these two times are independent of the shear rate $k_0$ within the shear band. However, this is a simple consequence of conservation of mass. The time $\hat{t}_c$ is the time at which enough small particles have dropped below $z=z_c$ to form a layer of small particles of depth $z_c$. These particles necessarily are transported from $z>z_c$, where their dynamics are independent of $k_0$. More precisely, 
the conservation law $\varphi_t+(a(z)f(\varphi))_z=0$ implies, together with the boundary conditions, that 
$
\int_0^1\varphi(z,t)\, dz
$
is independent of time. 
But $
\int_0^1\varphi(z,0)\, dz=1-z_0,
$
and for $t\geq\hat{t}_c$, we have $\varphi(z,t)=1, \, 0<z<z_c$, so that
$$
1-z_0=\int_0^1\varphi(z,t)\, dz=\int_0^{z_c}\varphi(z,t)\, dz+\int_{z_c}^{\gamma_1(t)}\varphi(z,t)\, dz=z_c+\int_{z_c}^{\gamma_1(t)}\varphi(z,t)\, dz.
$$
Consequently, $t=\hat{t}_c$ is the first time for which
$$
\int_{z_c}^{\gamma_1(t)}\varphi(z,t)\, dz=1-z_0-z_c,
$$
an equation that does not involve $k_0$, and which can be solved explicitly for $t=\hat{t}_c$ in the case $f(\varphi)=\varphi(\varphi-1)$.
Since $\hat{t}_c$ is independent of $k_0$, it follows that $t^*$ is as well.  

\subsection{ Case~II} 

\begin{figure}
\begin{center}
\includegraphics[width=.95\textwidth]{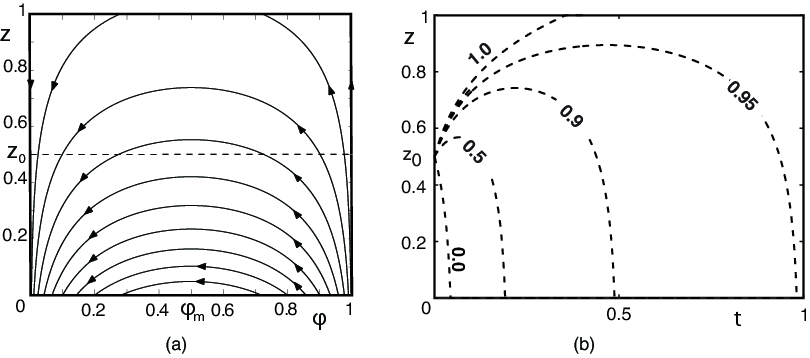}  
\caption{Characteristics (\ref{chars1_1}) in Case~II. (a) Phase portrait. \ (b) Projection onto $(t,z)$ plane. Curves are numbered by values of $\varphi_0$  in (\ref{zeqn}). } 
\label{phase_portrait} 
\end{center}
\end{figure}

First we examine the structure of the solution to (\ref{pde1}--\ref{bcs1}) for smooth functions $a(z)$ and fluxes $f(\varphi)$, satisfying the following conditions, consistent with the specifications of Case~II:
\beq\label{conditions}
f''(\varphi)>0, \quad f(0)=f(1)=0, \quad a(z)>0,\quad a'(z)<0.
\eeq
Under conditions (\ref{conditions}), the invariance (\ref{invariant}) of $a(z)f(\varphi)$ along characteristics is easily visualized, and gives the phase portrait for the vector field (\ref{chars1_1}), shown in Fig.~\ref{phase_portrait} (using $f(\varphi)=\varphi(\varphi-1),$ and parameter values (\ref{values})). 

From (\ref{conditions}), $f(\varphi)$ has a unique minimum, at  $\varphi_m\in (0,1)$. Trajectories of (\ref{chars1_1}) are horizontal at $\varphi_m$, and decreasing in $\varphi$, as shown in the figure. We also have $f'(\varphi)>0$ for $\varphi>\varphi_m$, so that (\ref{chars1_1}), (\ref{conditions}) imply that $z(t)$ is increasing there, and decreasing for $\varphi<\varphi_m$.

For fixed $z_0\in (0,1)$, the characteristic curves through $(z,t)=(z_0, 0)$ form a fan, with $\varphi=0$ corresponding to the line on the $z$ axis in the phase portrait of Fig.~\ref{phase_portrait}(a), and similarly $\varphi=1$ corresponding to the vertical line $\varphi=1$. As remarked earlier, these are the only straight characteristics, and the only characteristics on which $\varphi$ is constant. In between, the characteristics approach $z=0$ monotonically if $\varphi<\varphi_m$, or have a maximum in $z$ before reaching $z=0$, or reach $z=1$ monotonically, before the maximum is reached. This fan of characteristics forms the rarefaction wave emanating from $(z,t)=(z_0, 0)$, and joining $\varphi=0$ to $\varphi=1$. In Fig.~\ref{phase_portrait}(b), we show the characteristics in the $(z,t)$ plane, and remark that the pattern of characteristics is quite different from the pattern of contours of $\varphi$ in the representation in Fig.~\ref{solutions}(b). 

The rarefaction is reflected from the boundaries at $z=0, z=1$, forming a pair of shock waves \mbox{$z=\gamma_{\ell}(t), \ell=0,1$} that eventually meet at a time $t=t^*$, after which the solution consists of the single stationary shock from $\varphi=1$ to $\varphi=0$. Since $f(\varphi)$ is smooth and is zero at $\varphi=0,1$, let $f(\varphi)=\varphi(\varphi-1)g(\varphi)$, where $g(\varphi)$ is smooth and positive. Since $\varphi=1$ behind the shock $\gamma_0$, and $\varphi=0$ behind the shock $\gamma_1$, the ODE (\ref{rh1}) become
\beq
\frac{d\gamma_0}{dt}= a(\gamma_0) \varphi(\gamma_0,t) g(\varphi(\gamma_0,t)), \quad \frac{d\gamma_1}{dt}=a(\gamma_1) \left( \varphi(\gamma_1,t)-1 \right) g(\varphi(\gamma_1,t))
\eeq
In the general case, the values of $\varphi$ in the rarefaction are known only implicitly, and the shock waves can be determined only numerically. Even when these rarefaction values are given by a formula, as in the example presented below in connection with the experiment, the ODE may be intractable to explicit solution.

\subsubsection{Exponential Shear Rate} 

We are able to calculate the explicit solution, except for the reflected shocks, if we take $a(z)=a_0e^{- z/\lambda}$ and  $f(\varphi)=\varphi(\varphi-1)$. In Fig.~\ref{solutions}(b), we show the solution using parameter values (\ref{values}).

In the rarefaction wave, $\varphi$ is defined implicitly by (\ref{invariant}), which simplifies to the quadratic equation
\beq
\varphi(1-\varphi)=\varphi_o(1-\varphi_o)e^{(z-z_0)/\lambda }. \label{phiz}
\eeq
On the curves (\ref{phiz}), one for each $\varphi_o\in[0,1]$, $\varphi$ and $z$ evolve in time according to (\ref{chars1}), with initial conditions $\varphi(0)=\varphi_0$ and $z(0)=z_0$. In particular, the evolution of $\varphi$ is independent of the evolution of $z$, and in fact, $\varphi$ decays linearly along characteristics:
\beq
\frac{d\varphi}{dt}=\frac{a(z)}{\lambda}\varphi(\varphi-1)=-\frac{a_0}{\lambda}\varphi_o(1-\varphi_o)e^{-z_0/\lambda}
, \quad \varphi(0)=\varphi_0. \label{phit}
\eeq
Integrating yields 
\beq
\varphi(t)=-a'(z)\varphi_o(1-\varphi_o)e^{-z_0/\lambda}t+\varphi_o. \label{phit2}
\eeq
Substituting into equation (\ref{chars1}) 
gives the evolution of $z$:
\beq
\frac{dz}{dt}=a(z)(2\varphi-1)=a_0e^{-z/\lambda}(2\mathcal{P}t+2\varphi_o-1), \quad z(0)=z_0; \quad \mathcal{P}=-\frac{a_0}{\lambda} \varphi_o(1-\varphi_o)e^{-z_0/\lambda}. \label{zeqn}
\eeq
Therefore, 
\beq
e^{z/\lambda}-e^{z_0/\lambda}=\frac{a_0}{\lambda}(\mathcal{P}t^2+(2\varphi_o-1)t)=\frac{a_0}{\lambda}\left(-\frac{a_0}{\lambda}
\varphi_o(1-\varphi_o)e^{-z_0/\lambda}t^2+2\varphi_ot-t\right). \label{zeqn2}
\eeq
Eliminating $\varphi_0$ between (\ref{phit2}) and (\ref{zeqn2}) gives $\varphi=\varphi(z,t)$ in the rarefaction wave. Specifically, the quadratic equation (\ref{zeqn2}) can be solved for $\varphi_o=\varphi_o(z,t)$.  Then, $\varphi=\varphi(z,t)$ is obtained from (\ref{phit2}), with $\varphi_o=\varphi_o(z,t)$. The characteristics in Case~II using the experimentally determined parameters are shown in Fig.~\ref{phase_portrait}. 

Before specifying $\varphi(z,t)$ completely, we consider the leading edges of the wave emanating from $(z,t)=(z_0, 0)$. The leading edge $z = z_{min}(t)$ of the rarefaction approaching the lower boundary $z = 0$, carries $\varphi=\varphi_o=0$ and is the particle path of the first small particle to reach the lower boundary. Similarly, the leading edge $z = z_{max}(t)$ of the rarefaction, on which $\varphi=\varphi_o=1,$ approaches the upper boundary $z=1.$ Let $t_0, \, t_1$ be the times at which these curves reach $z=0, \, z=1$, respectively, so that $z_{min}(t_0)=0, \, z_{max}(t_1)=1$. These times are easily found from equation (\ref{zeqn2}) by substituting the relevant values for $z$ and $\varphi_o$:
\beq
t_0 = \frac{\lambda}{a_0}\left(e^{z_0/\lambda}-1\right), \qquad t_1 = \frac{\lambda}{a_0}\left(e^{1/\lambda} - e^{z_0/\lambda}\right).
\eeq
 
 Now we can solve the quadratic equation (\ref{zeqn2}) for $\varphi_o=\varphi_o(z,t)$, using the fact that $z= z_{max}(t) $ should give $\varphi_o=1$, to select the correct root of the equation:
 \beq
\varphi_o(z,t)=\frac{\lambda}{2a_0 t}\left(\frac{a_0 t}{\lambda} - 2e^{z_0/\lambda}+\sqrt{4e^{(z+z_0)/\lambda}+
\left(\frac{a_0 t}{\lambda}\right)^2}\right).
\label{phi01}    
\eeq
The entire rarefaction fan is now characterized using equation (\ref{phit2}) with $\varphi_o$ given by (\ref{phi01}):
\beq
\varphi(z,t)=-\frac{a_0}{\lambda} t \varphi_o(z,t)(1-\varphi_o(z,t))e^{-z_0/\lambda}+\varphi_o(z,t).
\label{phizt1}
\eeq

Next, we formulate an ODE for the reflected shocks $z=\gamma_0(t), \, z=\gamma_1(t)$,
originating from $(z,t)=(0, t_0)$, and from $(z,t)=(1, t_1)$ respectively. In order to track the shocks using the Rankine-Hugoniot condition, we use the expression for $\varphi=\varphi(z,t)$ given in equation (\ref{phizt1}) in the region between the shocks and between the outermost characteristics $ z=z_{min}(t)$ and $z=z_{max}(t)$, that is, the region in which we find the rarefaction fan. 

As in Case~I, the Rankine-Hugoniot condition for a shock curve $z=\gamma(t)$ for the PDE (\ref{pde1}) is given by (\ref{rh1}). 
Therefore, for $\gamma=\gamma_0(t)$, where $\varphi$ jumps from $\varphi=1$ to $\varphi=\varphi(\gamma_0,t)$, we have 
\beq
\frac{d\gamma_0}{dt}=a(\gamma_0)\varphi(\gamma_0,t), \quad t>t_0, \, \gamma_0(t_0)=0.
\label{gamma0}
\eeq
Similarly, for $\gamma=\gamma_1(t)$, where $\varphi$ jumps from $\varphi=\varphi(\gamma_1,t)$ to $\varphi=0$, we have 
\beq
\frac{d\gamma_1}{dt}=a(\gamma_1)(\varphi(\gamma_1,t)-1), \quad t>t_1, \, \gamma_1(t_1)=1.
\label{gamma1}
\eeq

We solve the equations (\ref{gamma0}) and (\ref{gamma1}) using the {\scshape Matlab} ODE solver  {\tt ode45}. The time $t^*$ at which $\gamma_0$ and $\gamma_1$ meet is calculated numerically.
 The entire solution  is shown in Fig.~\ref{solutions}(b), where the time scale has been normalized by the calculated $t^*=13.60;$ this is equivalent to setting the segregation rate $s=t^*$.

 \section{Discussion}\label{discuss} 

In both cases, the structure of the solutions in \S\ref{IBVP} captures the segregation process observed in the experiment, but there are significant differences between the two solutions, and between the solutions and the experiment.

Having computed solutions in each of Case~I and Case~II using the parameter values (\ref{values}) derived from experimental data,  and having noted the differences in construction, we are left with the observation that the time to segregation in the two cases is markedly different, with the ratio of final times $t^*$ being approximately $0.47$.  Because the exponential shear rate of Case~II is a better fit to the data, it is natural to regard it as the better model.

There are additional reasons to favor Case~II over Case~I. It would be tempting to blame the values of the shear rates $k_0$ and $k_1$ for the lack of agreement in $t^*$, but $t^*$ is in fact independent of $k_0$. Therefore, the problem has to lie in the upper part ($z>z_c$) of the sample. Indeed, the value of $t_1$ in Case~I is significantly smaller than in Case~II, and moreover, the shock $\gamma_1$ in Case~I is steeper over most of its path than the corresponding shock in Case~II. Both of these effects would be countered by decreasing the shear rate $k_1$ so that small particles approach the interface $z=z_c$ more slowly. However, the answer may be more subtle. In Case~I, $k_1 $ overestimates the experimental shear rate (and that of Case~II) near $z=1$, thereby promoting segregation there. This appears to be a more significant effect than the retarding of segregation closer to $z=z_c$ due to $k_1$ underestimating the shear rate there. It is also significant that, as indicated by the structure of characteristics in
Case~II (Fig.~\ref{phase_portrait}), the upper part of the domain strongly influences the segregation in the lower part, so that in Case~II, the entire domain is involved in determining the time to segregation. Consequently, the larger segregation rate near $z=0$ has significance in Case~II but none at all in Case~I. In summary, while the broad structure of the solutions (rarefaction wave and reflected shocks) is captured in both Case~I and Case~II, any quantitative 
comparison to experiments must utilize the more refined smooth fit to the shear rates of Case~II.

The comparison to experimental results faces a number of challenges, discussed in greater detail in our companion paper \cite{May-2009-SDP}. Here we note:
\begin{enumerate}

\item The model does not account for the rapid opening of void spaces due to Reynolds dilatency \cite{reynolds} during the initial transient dynamics. During this process, particles exhibit a more disorderly motion and the resulting measured velocity profile is inconsistent with those measured at later times. This short-time behavior undoubtedly accelerates the initial mixing of small and large particles prior to the establishment of the conditions described by the model. Therefore, the model only applies after this initial transient period, by which time the configuration of the particles is no longer given by (\ref{ic1}) and is instead a somewhat mixed state.

\item Full segregation is never achieved in the experiment, making it impossible to identify a final time  $t^*$ to compare with the models. To measure the progress of segregation in the experiment, we record the height of the top plate and relate the total volume of the sample to the degree of mixing/segregation. Because we observe that the volume approaches its final value only exponentially in time, we conjecture that isolated large particles remain trapped at the bottom and migrate upward on a slower timescale. This is a discrete effect not captured by a continuum model.

\item The model does not account for three-dimensional motion of the particles. While the assumption that the concentration of large and small particles depends only on depth is clearly unrealistic, $\varphi(z,t)$ might be considered reasonable as an average across horizontal layers of particles. The side walls of the experiment undoubtedly influence the dynamics, since they introduce lateral shear. It would be possible to include this effect in a multidimensional model, but verifying such a model with experiments would be hampered by the difficulty of tracking particles within the bulk.
\end{enumerate}

Finally, we note that the piecewise constant case (Case~I) has been studied both theoretically and numerically in various contexts, for example in flow in porous media composed of layers of different material such as sand and clay, and in connection with sedimentation  \cite{BachmannVovelle, Jimenez, KlingenbergRisebro, SeguinVovelle, WangSheng}. In these studies, a variety of techniques are introduced for studying entropy conditions and special solutions, as well as analysis of existence and uniqueness questions. Our application to segregation solves a special initial boundary value problem, but it reveals the unforeseen consequence that the material interface at $z=z_c$ removes the influence of the lower portion $z<z_c$ of the domain on the overall time scale of the evolution.

\subsection*{Acknowledgments}

The authors are grateful to Nico Gray for enlightening conversations about the model, and to Laura Golick and Katherine Phillips for assistance with the experiments. This research was supported by the National Science Foundation under grants  DMS-0604047 and DMS-0636590, and by NASA grant NNC04GB08G.


\end{document}